\documentclass[12pt]{article}%
\usepackage{graphicx}
\usepackage[intlimits]{amsmath}
\usepackage{latexsym}
\usepackage{amsfonts}
\usepackage{amssymb}%
 \makeatletter
\newfont{\footsc}{cmcsc10 at 8truept}
\newfont{\footbf}{cmbx10 at 8truept}
\newfont{\footrm}{cmr10 at 10truept}
\newtheorem{theorem}{Theorem}

\newtheorem{corollary}[theorem]{Corollary}

\newtheorem{definition}[theorem]{Definition}

\newenvironment{proof}[1][Proof]{\noindent{\textbf {#1}  }}  {\hfill$\Box$\bigskip}

\begin{document}

\title{The number of graphs with large forbidden subgraphs}
\author{B\'{e}la Bollob\'{a}s \thanks{Department of Pure Mathematics and Mathematical
Statistics, University of Cambridge, Cambridge CB3 0WB, UK and}
\thanks{Department of Mathematical Sciences, University of Memphis, Memphis TN
38152, USA} \thanks{Research supported in part by NSF grants DMS-0505550,
CNS-0721983, CCF-0728928 and DMS-0906634, and ARO grant W911NF-06-1-0076}
\ and Vladimir Nikiforov\footnotemark[2] \thanks{Research supported by NSF
Grant DMS-0906634.}}
\maketitle

\begin{abstract}
In this note, extending some results of Erd\H{o}s, Frankl, R\"{o}dl, Alexeev,
Bollob\'{a}s and Thomason, we determine asymptotically the number of graphs
which do not contain certain large subgraphs. In particular, we show that if
$H_{1},H_{2},\ldots$ are graphs with $\left\vert H_{n}\right\vert =o\left(
\log n\right)  $ and $\chi\left(  H_{n}\right)  =r_{n}+1,$ then the number
$S_{n}$ of graphs of order $n$ not containing $H_{n}$ satisfies
\[
\log_{2}S_{n}=\left(  1-1/r_{n}+o\left(  1\right)  \right)  \binom{n}{2}.
\]

We also give a similar statement for forbidden induced subgraphs.\medskip

\textbf{Keywords: }\textit{number of graphs; induced subgraph; removal lemma;
}

\end{abstract}

\subsection*{Introduction}

Given a graph $H,$ write $\mathcal{P}_{n}\left(  H\right)  $ for the set of
all labelled graphs of order $n$ not containing $H.$ In 1976, Erd\H{o}s,
Kleitman and Rothschild \cite{EKR76} gave a theorem implying that
\begin{equation}
\log_{2}\left\vert \mathcal{P}_{n}\left(  K_{r+1}\right)  \right\vert =\left(
1-1/r+o\left(  1\right)  \right)  \binom{n}{2}. \label{EKR}%
\end{equation}
In fact, this theorem, stated below, is considerably stronger than equation
(\ref{EKR}).\medskip

\textbf{Theorem A }\emph{Given }$r\geq2$\emph{ and }$\xi>0,$\emph{ there is
}$\rho=\rho\left(  r,\xi\right)  $\emph{ such that the number }$S_{n}$\emph{
of labelled graphs of sufficiently large order }$n$\emph{ containing at most
}$\rho n^{r+1}$\emph{ copies of }$K_{r+1}$\emph{ satisfies}%
\[
\left(  1-1/r\right)  \binom{n}{2}\leq\log_{2}S_{n}\leq\left(  1-1/r+\xi
\right)  \binom{n}{2}.
\]
\medskip

Ten years later, Erd\H{o}s, Frankl and R\"{o}dl \cite{EFR86} showed that the
conclusion in (\ref{EKR}) holds if $K_{r+1}$ is replaced by an arbitrary fixed
$\left(  r+1\right)  $-chromatic graph $H$. Surprisingly, Theorem A, combined
with an observation Erd\H{o}s \cite{Erd64} made in 1964 (Theorem E below),
easily imply the result of Erd\H{o}s, Frankl and R\"{o}dl, and even stronger
ones: for details see the first concluding remark at the end.

Here we shall give essentially best possible results that can be obtained by
replacing $H$ with a sequence of forbidden graphs whose order grows with $n$.
More precisely, given integers $r\geq2,$ $p\geq1,$ $q\geq1$ and real
$c\in\left(  0,1/2\right)  ,$ write $K_{r+1}\left(  p;q\right)  $ for the
complete $\left(  r+1\right)  $-partite graph with $r$ parts of size $p$ and
one part of size $q.$ Here and further, $\log$ with unspecified base stands
for the natural logarithm.

Our first result is the following theorem.

\begin{theorem}
\label{th1}Given $r\geq2$ and $0<\varepsilon\leq1/2,$ there is $\delta
=\delta\left(  \varepsilon\right)  >0$ such that for $n$ sufficiently large,%
\[
\left(  1-1/r\right)  \binom{n}{2}\leq\log_{2}\left\vert \mathcal{P}%
_{n}\left(  K_{r+1}\left(  \left\lfloor \delta\log n\right\rfloor ;\left\lceil
n^{1-\sqrt{\delta}}\right\rceil \right)  \right)  \right\vert \leq\left(
1-1/r+\varepsilon\right)  \binom{n}{2}.
\]

\end{theorem}

As we shall see, when $\varepsilon$ decreases, so does $\delta=\delta
(\varepsilon)$; this has the somewhat peculiar consequence that when
$\varepsilon$ decreases, the order of the forbidden graph $K_{r+1}\left(
\left\lfloor \delta\log n\right\rfloor ;\left\lceil n^{1-\sqrt{\delta}%
}\right\rceil \right)  $ increases; in fact, with the function $\delta
(\varepsilon)$ we shall take, this order is $\Theta\left(  n^{1-o\left(
\varepsilon\right)  }\right)  .$

Forgetting most of the vertices in the large vertex class of $K_{r+1}\left(
\left\lfloor \delta\log n\right\rfloor ;\left\lceil n^{1-\sqrt{\delta}%
}\right\rceil \right)  $, we get the following simplified assertion.

\begin{corollary}
Let $\left(  H_{n}\right)  $ be a sequence of graphs, with $\left\vert
H_{n}\right\vert =o\left(  \log n\right)  $ and $\chi\left(  H_{n}\right)
=r_{n}+1.$ Then, for every $\varepsilon>0$ and $n$ large enough,
\[
\left(  1-1/r_{n}\right)  \binom{n}{2}\leq\log_{2}\left\vert \mathcal{P}%
_{n}\left(  H_{n}\right)  \right\vert \leq\left(  1-1/r_{n}+\varepsilon
\right)  \binom{n}{2}.
\]

\end{corollary}

Indeed, if $r_{n}>1/\varepsilon,$ there is nothing to prove. If $r_{n}%
<1/\varepsilon$ then, as $H_{n}$ is a subgraph of the complete $\left(
r_{n}+1\right)  $-partite graph with all parts of size $\left\vert
H_{n}\right\vert =o\left(  \log n\right)  ,$ the upper bound follows when $n$
is sufficiently large. The lower bound follows as in Theorem \ref{th1}.

We should like to emphasize that Szemer\'{e}di's Regularity Lemma, a standard
tool to tackle questions like this, will not be used in our proof of Theorem
\ref{th1}, not even indirectly.

Next we turn to forbidden induced subgraphs, where the role of the chromatic
number is played by the \emph{coloring number }$\chi_{c}$, defined first in
\cite{BoTh95}, and given below.

\begin{definition}
Let $0\leq s\leq r$ be integers and let $\mathcal{H}\left(  r,s\right)  $ be
the class of graphs whose vertex sets can be partitioned into $s$ cliques and
$r-s$ independent sets. Given a graph property $\mathcal{P}$, the coloring
number $\chi_{c}\left(  \mathcal{P}\right)  $ is defined as%
\[
\chi_{c}\left(  H\right)  =\max\left\{  r:\mathcal{H}\left(  r,s\right)
\subset\mathcal{P}\text{ for some }s\in\left[  r\right]  \right\}
\]

\end{definition}

Also, given a graph $H,$ let us write $\mathcal{P}_{n}^{\ast}\left(  H\right)
$ for the set of graphs of order $n$ not containing $H$ as an induced
subgraph; clearly $\mathcal{P}_{n}^{\ast}\left(  H\right)  $ is a hereditary property.

A special case of a general result proved by Alexeev \cite{Ale93} and
independently by Bollob\'{a}s and Thomason \cite{BoTh95},\cite{BoTh97} is the
exact analogue of (\ref{EKR}): if $H$ is a fixed graph and $r=\chi_{c}\left(
\mathcal{P}_{n}^{\ast}\left(  H\right)  \right)  $, then
\begin{equation}
\log_{2}\left\vert \mathcal{P}_{n}^{\ast}\left(  H\right)  \right\vert
=\left(  1-1/r+o\left(  1\right)  \right)  \binom{n}{2}. \label{ABT}%
\end{equation}

Motivated by Theorem A, we first observe the following assertion, which is an
immediate consequence of the removal lemma of Alon, Fisher, Krivelevich and
Szegedy \cite{AFKS00} (Theorem B below) and the Alexeev-Bollob\'{a}s-Thomason
result (\ref{ABT}). We state it is as a theorem only to properly match Theorem
A for induced graphs.

\begin{theorem}
\label{th3}Let $H$ be a graph and let $r=\chi_{c}\left(  \mathcal{P}_{n}%
^{\ast}\left(  H\right)  \right)  .$ For every $\xi>0,$ there is a $\rho
=\rho\left(  H,\xi\right)  $ such that the number $S_{n}$ of graphs of
sufficiently large order $n$ containing at most $\rho n^{\left\vert
H\right\vert }$ induced copies of $H$ satisfies%
\[
\left(  1-1/r\right)  \binom{n}{2}\leq\log_{2}S_{n}\leq\left(  1-1/r+\xi
\right)  \binom{n}{2}.
\]

\end{theorem}

Note that our proof of this theorem uses implicitly Szemer\'{e}di's Regularity
Lemma. It would be interesting to find a proof avoiding this lemma. We know
from Erd\H{o}s, Kleitman and Rothschild \cite{EKR76} that this can be done
when $H=K_{r+1}.$\medskip

Next, we shall show that the conclusion in (\ref{ABT}) holds when $H$ is
replaced by a sequence of forbidden graphs whose order grows with $n$. To give
the precise statement, we need the following definition.

\begin{definition}
Given a labelled graph $H$ with $V\left(  H\right)  =\left[  h\right]  $ and
positive integers $p_{1},\ldots,p_{h}$, we say that a graph $F$ \emph{is of
type }$H\left(  p_{1},\ldots,p_{h}\right)  $ if $F$ can be obtained by
replacing each vertex $u\in V\left(  H\right)  $ with a graph $G_{u}$ of order
$p_{u}$ and each edge $uv\in E\left(  H\right)  $ with a complete bipartite
graph with vertex classes $V\left(  G_{u}\right)  $ and $V\left(
G_{v}\right)  $; if $uv\notin E\left(  H\right)  $ and $u\neq v,$ then $F$
contains no edges between $V\left(  G_{u}\right)  $ and $V\left(
G_{v}\right)  .$
\end{definition}

Now, given a labelled graph $H$ and positive integers $p$ and $q,$ let
\[
\mathcal{P}_{n}\left(  H;p,q\right)  =\left\{  G:%
\begin{array}
[c]{l}%
G\text{ is a labelled graph of order }n\text{ and }G\text{\ contains no
induced}\\
\text{subgraph of type }H\left(  p,\ldots,p,q\right)
\end{array}
\right\}  .
\]

Here is our second main result.

\begin{theorem}
\label{th2}Let $H$ be a labelled graph and let $r=\chi_{c}\left(
\mathcal{P}_{n}^{\ast}\left(  H\right)  \right)  .$ For every $\varepsilon>0,$
there is $\delta=\delta\left(  \varepsilon\right)  >0$ such that for $n$
sufficiently large%
\begin{equation}
\left(  1-1/r\right)  \binom{n}{2}\leq\log_{2}\left\vert \mathcal{P}%
_{n}\left(  H;\left\lfloor \delta\log n\right\rfloor ,\left\lceil
n^{1-\sqrt{\delta}}\right\rceil \right)  \right\vert \leq\left(
1-1/r+\varepsilon\right)  \binom{n}{2}. \label{bd1}%
\end{equation}

\end{theorem}

In some sense Theorems \ref{th1} and \ref{th2} are almost best possible, in
view of the following simple observation, that can be proved by considering
the random graph $G_{n,p}$ with $p\rightarrow1.$\medskip

\emph{Given }$r\geq2$\emph{ and }$\varepsilon>0,$\emph{ there is }$C>0$\emph{
such that the number of graphs }$S_{n}$\emph{ which do not contain }%
$K_{2}\left(  \left\lceil C\log n\right\rceil ,\left\lceil C\log n\right\rceil
\right)  $\emph{ satisfies }$S_{n}\geq\left(  1-\varepsilon\right)
2^{\binom{n}{2}}.$

\subsection*{Proofs}

For the proof of Theorem \ref{th3} we need a version of the Removal Lemma of
Ruzsa and Szemer\'{e}di \cite{RuSz75} for induced graphs; this result was
stated and proved in \cite{AFKS00}.\medskip

\textbf{Theorem B} \emph{Given a graph }$H$\emph{ and }$\alpha>0,$\emph{ there
is }$\beta=\beta\left(  \alpha\right)  >0$\emph{ such that if a graph }%
$G$\emph{ of order }$n$\emph{ contains fewer than }$\beta n^{\left\vert
H\right\vert }$\emph{ induced copies of }$H,$\emph{ then\ one can change at
most }$\alpha n^{2}$\emph{ edges of }$G$\emph{ so that the resulting graph
does not contain an induced copy of }$H.$\emph{\medskip}

We need also the following facts, which are Theorem 1 of \cite{Nik08a} and
Theorem 2 of \cite{Nik08b}.\medskip

\textbf{Theorem C} \emph{Let }$r\geq3,$\emph{ }$\left(  \ln n\right)
^{-1/r}\leq c\leq1/2,$\emph{ and let }$G$\emph{ be a graph with }$n$\emph{
vertices. If }$G$\emph{ contains more than }$cn^{r}$\emph{ copies of }$K_{r}%
$\emph{, then }$G$\emph{ contains a }$K_{r}\left(  s,\ldots s,t\right)
$\emph{ with }$s=\left\lfloor c^{r}\ln n\right\rfloor $\emph{ and
}$t>n^{1-c^{r-1}}.$\medskip

\textbf{Theorem D} \emph{Let }$2\leq h\leq n,$\emph{ }$\left(  \ln n\right)
^{-1/h^{2}}\leq c\leq1/4,$\emph{ let }$H$\emph{ be a graph of order }%
$h,$\emph{ and }$G$\emph{ be a graph of order }$n.$\emph{ If }$G$\emph{
contains more than }$cn^{h}$\emph{ induced copies of }$H,$\emph{ then }%
$G$\emph{ contains an induced subgraph of type }$H\left(  s,\ldots s,t\right)
,$\emph{ where }$s=\left\lfloor c^{h^{2}}\ln n\right\rfloor $\emph{ and
}$t>n^{1-c^{h-1}}.$\medskip

Note that Theorems \ref{th1}, \ref{th2}, and \ref{th3} have to be proved for
$n$ sufficiently large; thus, in the proofs below, we shall assume that $n$ is
as large as needed.\bigskip

\begin{proof}
[\textbf{Proof of Theorem \ref{th1}.}]Write $T_{r}\left(  n\right)  $ for the
$r$-partite Tur\'{a}n graph of order $n$ and note that no subgraph of
$T_{r}\left(  n\right)  $ contains a $K_{r+1}.$ Also, note that the number
$s_{n}^{\prime}$ of labelled spanning subgraphs of $T_{r}\left(  n\right)  $
satisfies
\[
\log_{2}s_{n}^{\prime}\geq\left(  1-1/r\right)  \frac{n^{2}}{2}-\frac{r}%
{8}\geq\left(  1-1/r\right)  \binom{n}{2},
\]
proving the lower bound in (\ref{bd1}); thus, to finish the proof of Theorem
\ref{th1} we need to prove the upper bound in (\ref{bd1}).

Fix $\varepsilon>0,$ let $\rho\left(  r,\cdot\right)  $ be the function from
Theorem A, and set $\delta=\rho\left(  r,\varepsilon\right)  ^{r+1}.$ If a
graph $G$ does not contain a $K_{r+1}\left(  \left\lfloor \delta\log
n\right\rfloor ;\left\lceil n^{1-\sqrt{\delta}}\right\rceil \right)  ,$ then
Theorem C implies that $G$ contains at most $\delta^{1/\left(  r+1\right)
}n^{r+1}=\rho\left(  r,\varepsilon\right)  n^{r+1}$ copies of $K_{r+1};$ in
turn, Theorem A implies that
\[
\log_{2}\left\vert \mathcal{P}_{n}\left(  K_{r+1}\left(  \left\lfloor
\delta\log n\right\rfloor ;\left\lceil n^{1-\sqrt{\delta}}\right\rceil
\right)  \right)  \right\vert \leq\left(  1-1/r+\varepsilon\right)  \binom
{n}{2},
\]
completing the proof of Theorem \ref{th1}.
\end{proof}

\begin{proof}
[\textbf{Proof of Theorem \ref{th3}.}]The lower bound is immediate since
$S_{n}$ must be at least as large as the number of graphs in $\mathcal{H}%
\left(  r,s\right)  $ of order $n$, and so, as in Theorem \ref{th1}, we see
that
\[
\log_{2}S_{n}\geq\left(  1-1/r\right)  \frac{n^{2}}{2}-\frac{r}{8}\geq\left(
1-1/r\right)  \binom{n}{2}.
\]

Let us now prove the upper bound. Fix $\varepsilon>0,$ and let $\sigma$ be
such that
\[
\frac{\varepsilon}{3}\geq\sigma\log_{2}\frac{4}{\sigma}.
\]
Let $\beta\left(  \cdot\right)  $ be the function of Theorem B, and set
$\delta=\beta\left(  \sigma/2\right)  .$ If a graph $G$ of order $n$ contains
at most $\delta n^{h}$ induced copies of $H,$ then Theorem B implies that all
induced copies of $H$ in $G$ can be destroyed by changing at most $\left(
\sigma/2\right)  n^{2}$ edges. Therefore, we see that
\begin{align*}
\log_{2}S_{n}  &  \leq\log_{2}\left\vert \mathcal{P}_{n}^{\ast}\left(
H\right)  \right\vert +\log_{2}\binom{\binom{n}{2}}{\sigma n^{2}/2}\leq
\log_{2}p_{n}+\sigma\frac{n^{2}}{2}\log_{2}\frac{4}{\sigma}\leq\log_{2}%
p_{n}+\frac{\varepsilon}{3}\frac{n^{2}}{2}\\
&  \leq\log_{2}p_{n}+\frac{\varepsilon}{2}\binom{n}{2}\leq\left(
1-1/r+\varepsilon\right)  \binom{n}{2}.
\end{align*}
The last inequality above follows from the Alexeev-Bollob\'{a}s-Thomason
result. This completes the proof of Theorem \ref{th3}.
\end{proof}

\begin{proof}
[\textbf{Proof of Theorem \ref{th2}.}]The lower bound is determined as in
Theorem \ref{th3}, so let us prove the upper bound. Let $\rho\left(
H,\cdot\right)  $ be the function from Theorem \ref{th3}. Fix $\varepsilon>0,$
let $\delta=\rho\left(  H,\varepsilon\right)  ^{h^{2}},$ and set
$p=\left\lfloor \delta\log n\right\rfloor ,$ $q=\left\lceil n^{1-\sqrt{\delta
}}\right\rceil .$ Suppose that a graph $G$ of order $n$ contains no induced
subgraph of type $H\left(  p,\ldots,p,q\right)  .$ Then, Theorem D implies
that $G$ contains at most $\delta^{1/h^{2}}n^{h}=\rho\left(  H,\varepsilon
\right)  n^{h}$ induced copies of $H.$ In turn, Theorem \ref{th3} implies that
the number $S_{n}$ of such graphs satisfies
\[
\log_{2}S_{n}\leq\left(  1-1/r+\varepsilon\right)  \binom{n}{2},
\]
completing the proof of Theorem \ref{th2}.
\end{proof}

\subsection*{Concluding remarks}

1. As mentioned at the beginning of this note, in \cite{Erd64}, equation
(18'), Erd\H{o}s gave a result about uniform hypergraphs, which implies the
following statement about graphs:\medskip

\textbf{Theorem E }\emph{Let }$r\geq2$. \emph{If a graph }$G$ \emph{of order
}$n$ \emph{contains at least }$\varepsilon n^{r}$\emph{ copies of }$K_{r},$
\emph{then }$G$ \emph{contains a copy of }%
\[
K_{r}\left(  \left\lfloor \delta\left(  \log n\right)  ^{1/\left(  r-1\right)
}\right\rfloor ,\ldots,\left\lfloor \delta\left(  \log n\right)  ^{1/\left(
r-1\right)  }\right\rfloor \right)
\]
\emph{for some }$\delta=\delta\left(  \varepsilon\right)  >0.$\medskip

In view of Theorem A, we immediately see the following corollary:\medskip

\emph{Given }$r\geq2$\emph{ and }$\varepsilon>0,$\emph{ there is }%
$\delta=\delta\left(  \varepsilon\right)  $\emph{ such that for }$n$\emph{
sufficiently large,}%
\begin{equation}
\log_{2}\left\vert \mathcal{P}^{n}\left(  K_{r+1}\left(  \left\lfloor
\delta\left(  \log n\right)  ^{1/r}\right\rfloor ,\ldots,\left\lfloor
\delta\left(  \log n\right)  ^{1/r}\right\rfloor \right)  \right)  \right\vert
\leq\left(  1-1/r+\varepsilon\right)  \binom{n}{2}. \label{E}%
\end{equation}
This statement could have been published as early as 1976, but the authors of
Theorem A somehow missed it, albeit Theorem E was used indeed in the proof of
Theorem A (see \cite{EKR76}, p. 20, line -5).\medskip

2. It is possible that the approach of \cite{Nik08a} can give an explicit
expression for $\delta\left(  \varepsilon\right)  $ in Theorem \ref{th1}. This
would help one to estimate how much $\left\vert \mathcal{P}_{n}\left(
K_{r+1}\left(  \left\lfloor \delta\log n\right\rfloor ;\left\lceil
n^{1-\sqrt{\delta}}\right\rceil \right)  \right)  \right\vert $ is larger than
the number of $r$-partite graphs of order $n$.\medskip

3. We reiterate the problem mentioned above: prove Theorem \ref{th3} avoiding
the use of Szemer\'{e}di's Regularity Lemma.\medskip

4. In the last two decades, the study of the number of graphs with given
properties has acquired a truly remarkable scale and sophistication, see,
e.g., \cite{ABBM09}\ and its references. Yet, we do not see a simple way to
accommodate the above results within this general framework.

\end{document}